\documentclass{amsart}
\usepackage[T1]{fontenc}
\usepackage{amssymb}
\usepackage{enumitem}
\usepackage{mathrsfs}
\usepackage{dsfont}
\usepackage[colorlinks, linkcolor=blue, citecolor=blue, urlcolor=blue]{hyperref}
\usepackage{tikz}

\newtheorem{theorem}{Theorem}
\newtheorem{lemma}[theorem]{Lemma}

\theoremstyle{definition}

\DeclareMathOperator{\ES}{ES}
\DeclareMathOperator{\ran}{ran}

\definecolor{diagramRegion}{HTML}{CEDFDA}
\definecolor{diagramFocus}{HTML}{9E8EB8}
\tikzset{
  proof diagram/.style={font=\small,line cap=round,line join=round},
  proof line/.style={black,line width=.55pt},
  proof axis/.style={black,line width=.45pt,->,>=stealth},
  proof guide/.style={black!68,line width=.4pt,densely dotted},
  proof cut/.style={black,line width=.45pt,dash pattern=on 2.5pt off 1.75pt}
}

\title{A new upper bound for the Erd\H{o}s--Szekeres conjecture}

\author{Yifan Hu}
\address[Yifan Hu]{School of Philosophy\\
Wuhan University\\
No.~299 Bayi Road\\
Wuhan\\
Hubei Province 430072\\
People's Republic of China}
\email{djfrankyifanhu@outlook.com}

\author{Guokuan Shao}
\address[Guokuan Shao]{School of Mathematics (Zhuhai)\\
Sun Yat-sen University\\
No.~2 Daxue Road\\
Zhuhai\\
Guangdong Province 519082\\
People's Republic of China}
\email{shaogk@mail.sysu.edu.cn}

\author{Guozhen Shen}
\address[Guozhen Shen]{Department of Philosophy (Zhuhai)\\
Sun Yat-sen University\\
No.~2 Daxue Road\\
Zhuhai\\
Guangdong Province 519082\\
People's Republic of China}
\email{shen\_guozhen@outlook.com}

\thanks{The central mathematical ideas are due to the authors.
The authors developed the technical approaches with the assistance of OpenAI's ChatGPT.
ChatGPT was also used for proofreading. The authors take full responsibility
for the paper's content and the accuracy of its references.
Guokuan Shao was partially supported by National Key R\&D Program of China grant number 2024YFA1015200
and National Natural Science Foundation of China grant number 12471082. Guozhen Shen was partially
supported by National Natural Science Foundation of China grant number 12671003.}



\begin{document}
\begin{abstract}
Let $\mathrm{ES}(n)$ denote the smallest integer such that every set of $\mathrm{ES}(n)$
points in the plane in general position contains $n$ points in convex position.
In 1935, Erd\H{o}s and Szekeres conjectured that $\mathrm{ES}(n)=2^{n-2}+1$,
and in 1960, they established the corresponding lower bound.
A breakthrough on the upper bound was achieved by Suk in 2017, who proved that
\[
\mathrm{ES}(n)\leqslant2^{n+O(n^{2/3}\log n)}.
\]
In this article, we further improve the upper bound by showing that
\[
\mathrm{ES}(n)\leqslant2^{n+O(n^{1/3})}.
\]
\end{abstract}
\maketitle

\section{Introduction}
For $n\geqslant3$, let $\mathrm{ES}(n)$ denote the smallest integer such that every set of $\mathrm{ES}(n)$
points in the plane in general position contains $n$ points in convex position.
Here, a set of points is said to be \emph{in general position} if no three of its points are collinear,
and \emph{in convex position} if every point is a vertex of its convex hull.
In their seminal 1935 paper~\cite{ES}, Erd\H{o}s and Szekeres conjectured that
\[
\ES(n)=2^{n-2}+1.
\]
They later established the corresponding lower bound by constructing a set of $2^{n-2}$
points in general position containing no $n$ points in convex position~\cite{ESlower}.
The problem of determining $\ES(n)$ is also known as the \emph{happy ending problem}.

After a suitable rotation, order the points from left to right.
A \emph{cap} is a chain with strictly decreasing edge slopes,
and a \emph{cup} is a chain with strictly increasing edge slopes.
The cup--cap recurrence of Erd\H{o}s and Szekeres~\cite{ES} yields
\[
\ES(n)\leqslant\binom{2n-4}{n-2}+1.
\]
Using a projective transformation, T\'oth and Valtr~\cite{TV} reduced this bound by approximately a factor of two.

Throughout this article, all logarithms are to base two. In 2017, Suk~\cite{Suk} obtained the remarkable upper bound
\[
\mathrm{ES}(n)\leqslant2^{n+O(n^{2/3}\log n)}.
\]
He used the positive fraction theorem of P\'or and Valtr~\cite{PV} to find a cap or cup such that
the exterior support region associated with each edge contains a fixed fraction of the original point set,
with the fraction depending only on the length of the chain.
He then applied Dilworth's theorem within each region and joined convex chains from different regions.
By counting compatible chains in adjacent or alternating regions,
Holmsen, Mojarrad, Pach, and Tardos~\cite{HMPT} improved the upper bound to
\[
\mathrm{ES}(n)\leqslant2^{n+O(\sqrt{n\log n})}.
\]
In their basic estimate with $k$ regions, the error is $8k+\frac{2n\log n}{k}$, whose minimum is of order $\sqrt{n\log n}$.

The upper and lower hulls of a convex $n$-gon form a cap and a cup, respectively, with common endpoints and a total of $n$ edges.
Using a lattice-path counting argument, Baek and Balko~\cite{BB} proved that $2^{n-2}+1$
points suffice when only the right endpoints of the two chains are required to coincide.
For a nonempty finite set of points $P$, let $\chi(P)$ denote the maximum size of a subset of $P$ in convex position.
In this article, we establish the following bound.

\begin{theorem}\label{main}
There exists an absolute constant $C>0$ such that, for every nonempty finite set of points $P$ in the plane in general position,
\[
\log|P|\leqslant\chi(P)+C\chi(P)^{1/3}.
\]
Consequently, $\ES(n)\leqslant 2^{n+O(n^{1/3})}$ for all integers $n\geqslant3$.
\end{theorem}

The novelty of our proof lies in combining a geometric reduction with an adaptive selection of lattice paths.
For a nonempty subset $X$ of $P$, let $S(X)$ be the sum of the maximum lengths, measured in edges, of a cap and a cup in $X$.
We encode each point of $X$ by a distinct lattice path of length $S(X)$, thereby translating
geometric information about $X$ into combinatorial information about its associated paths.
If these paths are sufficiently close to one another, we can join caps and cups from the left and right
halves of $X$ along supporting edges, while controlling the number of vertices lost in the process.
Since the resulting cap and cup together form a convex polygon, we can find a large subset
of one of the two halves for which the excess $S(X)-\chi(P)$ is substantially smaller.

The main difficulty is to select a large family of paths satisfying the condition required by the geometric reduction.
Our idea is to describe the paths by recording their heights at suitably chosen positions.
Rather than choosing these positions in advance, we construct a common partition adaptively,
subdividing an interval only when the heights of most of the remaining paths vary considerably on it.
An elementary counting argument shows that relatively few paths exhibit large variations on many disjoint intervals,
thereby limiting the number of subdivisions needed. Consequently, by recording only a small amount of information about each path,
we can select a large subfamily of paths that are sufficiently close to one another while keeping the loss in logarithmic size small.
Starting with $X=P$, we repeatedly perform this selection and apply the geometric reduction,
adjusting the accuracy of the path selection to the current excess $S(X)-\chi(P)$.
The key advantage of this approach is that the excess decreases geometrically, while the total loss
in logarithmic size throughout the iteration is only $O(\chi(P)^{1/3})$, yielding the desired bound.

We conclude the introduction by mentioning several extensions of the happy ending problem.
Fox, Pach, Sudakov, and Suk~\cite{FPSS} studied extensions
to monotone paths in ordered uniform hypergraphs and to families of convex bodies.
Conlon, Fox, He, Mubayi, Suk, and Verstra\"ete~\cite{CFHMSV}
considered point sets not necessarily in general position,
while Pohoata and Zakharov \cite{PZ} obtained subexponential
bounds for the happy ending problem in higher dimensions.

\section{The geometric reduction theorem}
Fix a nonempty finite set $P$ of points in the plane in general position.
Choose a coordinate system in which the points of $P$ have distinct first coordinates.
We write $p\prec q$ if the first coordinate of $p$ is smaller than that of $q$,
and let $\sigma(p,q)$ denote the slope of the line through $p$ and $q$.
A chain $p_0\prec\dots\prec p_k$ has \emph{length} $k$.
It is a \emph{cap} if its successive edge slopes are strictly decreasing,
and a \emph{cup} if they are strictly increasing.
Chains of length zero or one are both caps and cups.

For a nonempty set $X\subseteq P$, let $A(X)$ and $U(X)$ denote the maximum lengths of a cap and a cup in $X$,
respectively, and put $S(X)=A(X)+U(X)$. The union of a cap and a cup with common distinct endpoints is in convex position.
Unless both chains consist of the same edge, they form the upper and lower boundaries of a convex polygon, respectively.
Their vertex sets intersect only at their common endpoints, so the number of vertices in their union equals the sum of their lengths.
A cap or cup of length $k$ itself consists of $k+1$ vertices in convex position. Hence
\[
A(X),U(X)\leqslant\chi(X)-1,\qquad
S(X)\leqslant2\chi(X)-2.
\]

Fix a nonempty set $X\subseteq P$. For $p\in X$ and $\theta\in\mathds R$, let $x_p(\theta)$ denote
the maximum length of a cap in $X$ ending at $p$ whose last edge has slope at least $\theta$.
Similarly, let $y_p(\theta)$ denote the maximum length of a cup in $X$ ending at $p$
whose last edge has slope at most $\theta$, and put $z_p(\theta)=x_p(\theta)+y_p(\theta)$.
Chains of length zero are regarded as satisfying both slope conditions.
Note that these ranks depend on $X$, although this dependence is suppressed in the notation.

A lattice path consists of unit east and north steps.
We define the \emph{time} and \emph{height} of a vertex $(i,j)$ to be $i+j$ and $i-j$, respectively.
Vertices are compared with respect to the coordinatewise order; two vertices are \emph{incomparable}
if neither is less than or equal to the other.

For each $p\in X$, as $\theta$ decreases, both coordinates of $(x_p(\theta),U(X)-y_p(\theta))$ are nondecreasing.
Extend the resulting finite chain of vertices to a lattice path from $(0,0)$ to $(A(X),U(X))$.
We call this assignment the \emph{ending-rank encoding}; related encodings appear in~\cite{MS,NY}.

For any $p,q\in X$ with $p\prec q$, no edge ending at $p$ has slope $\sigma(p,q)$, by general position.
Appending $q$ to a cap or cup attaining the corresponding rank at $p$ therefore gives
\begin{equation}\label{eq:rank-increments}
x_q(\sigma(p,q))\geqslant x_p(\sigma(p,q))+1,\qquad
y_q(\sigma(p,q))\geqslant y_p(\sigma(p,q))+1.
\end{equation}
Thus the rank vertex corresponding to $q$ has a larger first coordinate
and a smaller second coordinate than that corresponding to $p$.
These two vertices are incomparable and hence cannot lie on the same lattice path.
Their heights are $z_p(\sigma(p,q))-U(X)$ and $z_q(\sigma(p,q))-U(X)$, respectively,
and so differ by $z_q(\sigma(p,q))-z_p(\sigma(p,q))$.
In particular, distinct points of $X$ are assigned distinct lattice paths.
Since there are exactly $\binom{A(X)+U(X)}{A(X)}$ lattice paths with the prescribed endpoints, it follows that
\begin{equation}\label{eq:encoding-bound}
|X|\leqslant\binom{A(X)+U(X)}{A(X)}\leqslant2^{S(X)}.
\end{equation}

\begin{theorem}[geometric reduction theorem]\label{prop:reduction}
Let $d\geqslant1$ be an integer and let $T$ be a nonempty subset of $X$. Suppose that
\begin{equation}\label{eq:rank-width}
z_q(\sigma(p,q))-z_p(\sigma(p,q))\leqslant3d
\qquad\text{for all }p,q\in T\text{ with }p\prec q.
\end{equation}
Then there exists a nonempty set $X'\subseteq T$ such that
\begin{equation}\label{eq:reduction}
|X'|\geqslant\frac{|T|}{3\lceil\chi(X)/d\rceil},\qquad
S(X')\leqslant\frac{\chi(X)+S(X)}2+4d.
\end{equation}
\end{theorem}
\begin{proof}
The case $|T|=1$ is immediate. Otherwise, partition $T$ into the set $T_-$
of its first $\lfloor|T|/2\rfloor$ points and the set $T_+$ of its remaining points.
Each of these sets contains at least $|T|/3$ points. Let $k_-$ denote the maximum sum
of the lengths of a cap and a cup in $T_-$ with a common starting point. Similarly,
let $k_+$ denote the maximum sum of the lengths of a cap and a cup in $T_+$
with a common ending point. We claim that
\begin{equation}\label{eq:half-sums}
k_-+k_+\leqslant\chi(X)+S(X)+6d.
\end{equation}

Choose pairs of chains attaining these maxima, with common starting point $p$
on the left and common ending point $q$ on the right.
Consider the upper convex hull of the union of the vertex sets of the two caps.
Let $ef$ be its edge crossing the gap between $T_-$ and $T_+$, with $e\in T_-$ and $f\in T_+$.
The line through $e$ and $f$ is an upper supporting line of both caps,
and every other vertex lies strictly below it by general position. Note that $e=p$ or $f=q$ may occur.
Put $\theta=\sigma(e,f)$. Concatenate the portion of the left cap from $p$ to $e$, the edge $ef$,
and the portion of the right cap from $f$ to $q$. The slope of the edge entering $e$ is greater than $\theta$,
and that of the edge leaving $f$ is smaller than $\theta$, whenever these edges exist.
Thus the resulting chain is a cap; see Figure~\ref{fig:supporting-edge}.

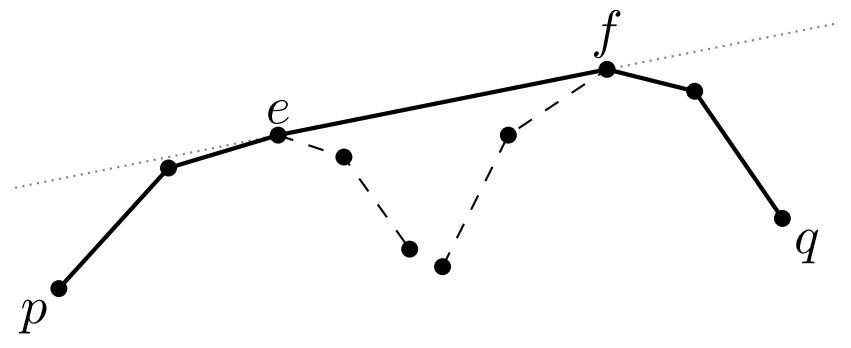
\begin{figure}[h]
\centering
\begin{tikzpicture}[scale=1.45]
 \coordinate (p) at (0,0);
 \coordinate (lone) at (0.5,0.55);
 \coordinate (e) at (1,0.7);
 \coordinate (lthree) at (1.3,0.6);
 \coordinate (lfour) at (1.6,0.18);
 \coordinate (rzero) at (1.75,0.1);
 \coordinate (rone) at (2.05,0.7);
 \coordinate (f) at (2.5,1);
 \coordinate (rthree) at (2.9,0.9);
 \coordinate (q) at (3.3,0.32);
 \draw[densely dotted,gray] (-0.2,0.46)--(3.55,1.21);
 \draw[thick] (p)--(lone)--(e)--(f)--(rthree)--(q);
 \draw[dashed] (e)--(lthree)--(lfour);
 \draw[dashed] (rzero)--(rone)--(f);
 \foreach \vertex in {p,lone,e,lthree,lfour,rzero,rone,f,rthree,q}
   \fill (\vertex) circle (1.1pt);
 \node[below left,inner sep=2pt] at (p) {$p$};
 \node[above,inner sep=2pt] at (e) {$e$};
 \node[above,inner sep=2pt] at (f) {$f$};
 \node[below right,inner sep=2pt] at (q) {$q$};
\end{tikzpicture}
\caption{The upper supporting edge joins the solid portions of the two caps.
The discarded portions have lengths at most $A(X)-x_e(\theta)$ and $x_f(\theta)$, respectively.}
\label{fig:supporting-edge}
\end{figure}

The discarded left suffix, if nontrivial, starts at $e$ with an edge of slope less than~$\theta$.
It can therefore be appended to a cap attaining $x_e(\theta)$, so its length is at most $A(X)-x_e(\theta)$.
The discarded right prefix, if nontrivial, ends at $f$ with an edge of slope greater than $\theta$,
so its length is at most $x_f(\theta)$. The total number of discarded cap edges is therefore at most
\[
A(X)+x_f(\theta)-x_e(\theta)\leqslant A(X)+z_f(\theta)-z_e(\theta)\leqslant A(X)+3d,
\]
where the inequalities follow from \eqref{eq:rank-increments} and \eqref{eq:rank-width}.
Reflection across a horizontal line gives the analogous bound $U(X)+3d$ for the number of discarded cup edges.
The resulting cap and cup have common endpoints $p,q$, so the sum of their lengths is at most $\chi(X)$.
Since the joining procedure adds two edges, we obtain
\[
k_-+k_+\leqslant\chi(X)+S(X)+6d-2.
\]
This proves \eqref{eq:half-sums}.

Choose a half $H$ whose corresponding maximum $k$ satisfies
\[
k\leqslant\frac{\chi(X)+S(X)}2+3d.
\]
For $p\in H$, let $a(p)$ and $u(p)$ denote the maximum lengths of a cap and a cup in $H$,
respectively, starting at $p$ if $H=T_-$, and ending at $p$ if $H=T_+$.
Then $a(p)+u(p)\leqslant k$ and $0\leqslant u(p)\leqslant\chi(X)-1$.
Partition $H$ into classes according to the value of $\lfloor u(p)/d\rfloor$,
and let $X'$ be a largest class, with common value $r$.
Since there are at most $\lceil\chi(X)/d\rceil$ classes, we have
\[
|X'|\geqslant\frac{|H|}{\lceil\chi(X)/d\rceil}\geqslant\frac{|T|}{3\lceil\chi(X)/d\rceil}.
\]
For any cap in $X'$, let $p$ be its starting point if $H=T_-$ and its ending point if $H=T_+$.
Its length is at most $a(p)\leqslant k-u(p)\leqslant k-rd$.
Similarly, every cup in $X'$ has length at most $(r+1)d-1$. Thus
\[
S(X')\leqslant k+d-1\leqslant\frac{\chi(X)+S(X)}{2}+4d,
\]
which completes the proof.
\end{proof}

\section{The path selection theorem}
In this section, we study combinatorial questions concerning lattice paths.
For a sign sequence $w=(w_1,\dots,w_n)\in\{-1,1\}^n$, associate an east step with $+1$ and a north step with $-1$, and put
\[
H_w(0)=0,\qquad H_w(t)=\sum_{i=1}^t w_i\qquad(1\leqslant t\leqslant n).
\]
The vertex at time $t$ of the corresponding lattice path starting at $(0,0)$ has height $H_w(t)$ and coordinates
\[
\left(\frac{t+H_w(t)}2,\frac{t-H_w(t)}2\right).
\]
Clearly, two vertices at times $s,t$ and of heights $g,h$, respectively, are incomparable if and only if $|g-h|>|s-t|$.

A family of lattice paths has \emph{width at most $d$} if any two incomparable vertices
on paths in the family differ in height by at most $d$. Equivalently,
\[
|H_w(s)-H_{w'}(t)|\leqslant\max\{d,|s-t|\}
\]
for all $w,w'$ in the family and all $0\leqslant s,t\leqslant n$.
Throughout this section, all time variables and interval endpoints are integers.
For an interval $I=[s,t]$ with $0\leqslant s\leqslant t\leqslant n$, we define the \emph{range} of $H_w$ on $I$ by
\[
\ran_w(I)=\max_{s\leqslant u\leqslant t}H_w(u)-\min_{s\leqslant u\leqslant t}H_w(u).
\]

\begin{lemma}\label{lem:excursion-count}
Let $d,m\geqslant1$ be powers of two with $m\leqslant d^2/64$.
There are at most $2^{m-d^2/(32m)}$ sequences $w\in\{-1,1\}^m$ satisfying $\ran_w([0,m])>d$.
\end{lemma}
\begin{proof}
If $m\leqslant d$, then $\ran_w([0,m])\leqslant m\leqslant d$ for every $w\in\{-1,1\}^m$,
and the assertion is immediate. Assume that $m>d$, and put
\[
b=\frac{8m}{d},\qquad k=\frac{d^2}{64m},\qquad h=\frac d2.
\]
By the hypotheses, $b$ is an even integer and $k,h$ are positive integers. Moreover,
\[
m=kb^2,\qquad h=4kb.
\]

We shall use the elementary inequalities
\begin{equation}\label{eq:colour-elementary}
(b+1)^b\geqslant2b^b,\qquad
(b^2-1)^{b^2/2}\geqslant\tfrac12(b^2)^{b^2/2}.
\end{equation}
For the first inequality, among words of length $b$ over an alphabet of $b+1$ letters,
count those in which a specified letter occurs at most once. Their number is $b^b+bb^{b-1}=2b^b$.
For the second inequality, put $r=b^2/2$. Among words of length $r$ over an alphabet of $b^2$ letters,
at most $r(b^2)^{r-1}=(b^2)^r/2$ contain a specified letter, as can be seen by summing over its possible positions.
Thus at least half of the words avoid that letter, which proves the second inequality.

Let $\mathcal H$ be the family of sequences $w\in\{-1,1\}^m$ for which $H_w(t)=h$ for some $t$.
We bound $|\mathcal H|$ by counting coloured sign sequences.

For each $w\in\{-1,1\}^m$, colour its entries according to the following rule.
Until and including the first time at which $H_w$ reaches $h$,
each positive entry has $b+1$ available colours and each negative entry has $b-1$ available colours.
All subsequent entries have $b$ available colours, regardless of their signs.
If $H_w$ never reaches $h$, the first rule applies to all entries.

Given any admissible coloured prefix of length less than $m$,
there are exactly $2b$ choices for the next sign together with its colour:
before the first visit to $h$, there are $(b+1)+(b-1)$ choices,
and afterwards there are $b+b$ choices.
Hence, there are exactly $(2b)^m$ coloured sequences in total.

Fix $w\in\mathcal H$. Suppose that the prefix ending at the first time at which $H_w$ reaches $h$
contains $i$ negative entries. It then contains $h+i$ positive entries. Hence, the number of colourings of $w$ is
\[
(b+1)^{h+i}(b-1)^ib^{m-h-2i}=b^m(1+1/b)^h(1-1/b^2)^i\geqslant b^m2^{4k}2^{-k}=b^m2^{3k},
\]
where the inequality follows from~\eqref{eq:colour-elementary}, together with $h=4kb$ and $i\leqslant m/2=kb^2/2$.

Since distinct sign sequences have disjoint sets of colourings, we obtain
\[
|\mathcal H|b^m2^{3k}\leqslant(2b)^m,
\]
and consequently
\[
|\mathcal H|\leqslant2^{m-3k}.
\]
By symmetry, the same bound holds for the number of sequences whose height reaches $-h$.

If $\operatorname{ran}_w([0,m])>d=2h$, then some height is greater than $h$ or smaller than $-h$.
Since $H_w(0)=0$ and successive heights differ by $1$, $H_w$ must reach either $h$ or $-h$.
Therefore, the number of such sequences is at most
\[
2\cdot2^{m-3k}\leqslant2^{m-2k}=2^{m-d^2/(32m)},
\]
as required.
\end{proof}

\begin{theorem}[path selection theorem]\label{thm:paths}
There exists an absolute constant $C>0$ with the following property.
Let $d\geqslant1$ be a power of two, let $n$ be a positive integer with $n\leqslant2d^3$,
and suppose that $\Omega$ is a nonempty subset of $\{-1,1\}^n$ satisfying $|\Omega|\geqslant2^{n-32d}$.
Then there exists a nonempty subfamily $\Gamma\subseteq\Omega$ of width at most $3d$ such that
\begin{equation}\label{eq:path-loss}
\log\frac{|\Omega|}{|\Gamma|}\leqslant C\left(\sqrt{\frac nd}+\log(d+1)\right).
\end{equation}
\end{theorem}
\begin{proof}
If $n\leqslant d$, then the whole family has width at most $2n\leqslant2d$. Thus we may assume that $n>d$.
Let $L$ be the largest power of two not exceeding $\max\{1,\sqrt{nd}/256\}$, and put $R=\lceil n/L\rceil$ and $N=RL$. Then
\begin{equation}\label{eq:scale}
\begin{gathered}
\frac{\sqrt{nd}}{512}\leqslant L
\leqslant\max\left\{1,\frac{\sqrt{nd}}{256}\right\},\\
R=O\left(\sqrt{\frac nd}\right)=O(d),\qquad
\frac Ld=O\left(\sqrt{\frac nd}\right).
\end{gathered}
\end{equation}
For each $w\in\Omega$, extend $H_w$ constantly from $n$ to $N$. This preserves the inequality
\[
|H_w(s)-H_w(t)|\leqslant|s-t|
\]
for all $0\leqslant s,t\leqslant N$. The \emph{length} of an interval $I=[s,t]$ is defined to be $t-s$.

We first construct a partition of $[0,N]$ together with a subfamily of $\Omega$
whose members have range at most $d$ on each interval of the partition.
Start with the partition into $R$ intervals of length $L$ and the family $\Omega$.
Examine the intervals from left to right. If at least half of the current family consists of sequences $w$
satisfying $\operatorname{ran}_w(I)\leqslant d$ on the interval $I$ under consideration,
retain these sequences and keep~$I$. Otherwise, retain the sequences $w$ satisfying $\operatorname{ran}_w(I)>d$,
bisect $I$, and examine the resulting two intervals. Since an interval of length at most $d$ is never bisected,
the procedure terminates with a nonempty family $\Omega'$ and a partition $\mathcal I$ of $[0,N]$ into, say, $q$ intervals.

Every $w\in\Omega'$ has range at most $d$ on each interval of $\mathcal I$ and range greater than $d$
on each bisected interval. Since there are $q-R$ bisections and $q$ final intervals, the total number
of examinations is $2q-R$. Each examination retains at least half of the current family, so
\begin{equation}\label{eq:survivors}
|\Omega'|\geqslant|\Omega|2^{-(2q-R)}.
\end{equation}

We next bound the number of bisections. Put $b=q-R$ and $D=n-\log|\Omega'|$. By the hypothesis on $|\Omega|$ and~\eqref{eq:survivors},
\begin{equation}\label{eq:deficit}
D\leqslant32d+R+2b.
\end{equation}
If $b>0$, then $L>d$, and hence
\begin{equation}\label{eq:small-L}
L\leqslant\frac{\sqrt{nd}}{256}\leqslant\frac{\sqrt2d^2}{256}<\frac{d^2}{128}.
\end{equation}

For each power of two $m\leqslant L$, let $b_m$ denote the number of bisected intervals of length $m$.
These intervals have pairwise disjoint interiors, and every $w\in\Omega'$ has range greater than $d$ on each of them.
For every such interval $I$, the range of the constant extension of $H_w$ on $I$ is already attained on $I\cap[0,n]$.
Indeed, $I$ cannot lie wholly to the right of $n$, since its range would then be zero; and if $I$ crosses $n$,
the constant part only repeats the value $H_w(n)$.

Extend each sequence in $\Omega'$ by a sign sequence of length $N-n$ in all possible ways.
By the preceding observation, each of the resulting $|\Omega'|2^{N-n}$ sequences
still has range greater than $d$ on each of the $b_m$ intervals. If $b_m>0$,
then \eqref{eq:small-L} ensures that Lemma~\ref{lem:excursion-count} applies.
The range on an interval $[s,t]$ depends only on the increments $w_{s+1},\ldots,w_t$.
Since the intervals have pairwise disjoint interiors, they impose conditions on disjoint sets of entries, and hence
\[
|\Omega'|2^{N-n}\leqslant 2^{N-d^2b_m/(32m)}.
\]
Consequently,
\begin{equation}\label{eq:bisection-count}
|\Omega'|\leqslant2^{n-d^2b_m/(32m)},\qquad
b_m\leqslant\frac{32mD}{d^2},
\end{equation}
which also hold trivially when $b_m=0$. Summing over all possible lengths $m=L,L/2,\dots,1$, we obtain
\begin{equation}\label{eq:total-bisections}
b=\sum_m b_m\leqslant\frac{32D}{d^2}\sum_m m\leqslant\frac{64LD}{d^2}.
\end{equation}
If $b>0$, combining \eqref{eq:deficit}, \eqref{eq:small-L}, and \eqref{eq:total-bisections} gives
\[
D\leqslant32d+R+\frac{128LD}{d^2}\leqslant32d+R+\frac{D}{\sqrt2}.
\]
Since $R=O(d)$, it follows that $D=O(d)$. If $b=0$, the same conclusion follows directly from~\eqref{eq:deficit}.
Substituting this bound into~\eqref{eq:total-bisections}, we obtain
\begin{equation}\label{eq:interval-count}
q=R+b=O\left(R+\frac Ld\right)=O\left(\sqrt{\frac nd}\right).
\end{equation}

We now select a large subfamily of $\Omega'$ by recording the heights
of the paths at the left endpoints of the intervals in $\mathcal I$.
Call an interval in $\mathcal I$ \emph{short} if its length is less than $2d$, and \emph{long} otherwise.
We claim that there are only $O(1)$ short intervals. Indeed, if $L\leqslant d$, then no interval is bisected,
and~\eqref{eq:scale} implies that $\sqrt{n/d}\leqslant512$. Thus $q=R=O(1)$.
Otherwise, since all interval lengths are powers of two and no interval of length at most $d$ is bisected,
every short interval has length $d$ and arises from the bisection of an interval of length $2d$.
By~\eqref{eq:bisection-count}, the number of short intervals is therefore at most
\[
2b_{2d}\leqslant\frac{128D}{d}=O(1).
\]

We record the initial heights exactly on short intervals,
so that paths with the same records agree at their left endpoints.
On long intervals, we round the initial heights down to
multiples of $d$ to reduce the number of possible records.
For each $w\in\Omega'$ and $I=[s,t]\in\mathcal I$, define
\[
\lambda_w(I)=
\begin{cases}
H_w(s) & \text{if }I\text{ is short,}\\
d\lfloor H_w(s)/d\rfloor & \text{if }I\text{ is long.}
\end{cases}
\]
In either case,
\begin{equation}\label{eq:hd}
0\leqslant H_w(s)-\lambda_w(I)<d.
\end{equation}
We refer to $\lambda_w(I)$ as the \emph{record} of $w$ on $I$.
Since every $w\in\Omega'$ has range at most $d$ on each interval of $\mathcal I$,
the records on consecutive intervals differ by at most $2d$.
Thus, once the record on the preceding interval is fixed,
there are at most five possible records on a long interval and at most $4d+1$ possible records on a short interval.
The record on the first interval is zero. Since there are $q$ intervals and only $O(1)$ short intervals,
the number of distinct record sequences is at most
\[
5^q(4d+1)^{O(1)}.
\]
Let $\Gamma$ be a largest subfamily of $\Omega'$ whose members have
identical records on every interval of $\mathcal I$. By~\eqref{eq:survivors} and~\eqref{eq:interval-count},
\[
\begin{aligned}
\log\frac{|\Omega|}{|\Gamma|}
& =\log\frac{|\Omega|}{|\Omega'|}+\log\frac{|\Omega'|}{|\Gamma|}\\
& \leqslant2q-R+q\log5+O(\log(d+1))\\
& =O\left(\sqrt{\frac nd}+\log(d+1)\right).
\end{aligned}
\]

It remains to show that $\Gamma$ has width at most $3d$. By~\eqref{eq:hd} and the choice of $\Gamma$,
the heights of any two paths in $\Gamma$ differ by at most $d$ at the left endpoint of each interval in $\mathcal I$,
and agree if the interval is short.

Fix $w,w'\in\Gamma$ and $0\leqslant s\leqslant t\leqslant n$.
Suppose first that $[s,t]$ contains an entire interval $I=[u,v]\in\mathcal I$.
If $I$ is short, then $H_w(u)=H_{w'}(u)$, so
\[
|H_w(s)-H_{w'}(t)|\leqslant(u-s)+(t-u)=t-s.
\]
If $I$ is long, then the range bound on $I$ gives
\[
|H_w(s)-H_{w'}(t)|\leqslant(u-s)+2d+(t-v)\leqslant t-s,
\]
since $v-u\geqslant2d$.

Otherwise, $s,t$ lie in the same interval or in two adjacent intervals.
Let $u$ be the left endpoint of that interval in the former case and
the common endpoint of the two intervals in the latter case.
By the range bounds and the agreement of records,
\[
\begin{aligned}
|H_w(s)-H_{w'}(t)|
& \leqslant|H_w(s)-H_w(u)|+|H_w(u)-H_{w'}(u)|+|H_{w'}(u)-H_{w'}(t)|\\
& \leqslant3d.
\end{aligned}
\]

In either case, we have
\[
|H_w(s)-H_{w'}(t)|\leqslant\max\{3d,t-s\}.
\]
Thus $\Gamma$ has width at most $3d$, as required.
\end{proof}

\section{Proof of the main theorem}
We are now ready to prove our main theorem.
\begin{proof}[Proof of Theorem~\ref{main}]
We proceed by constructing a decreasing sequence $(X_i)$ of nonempty subsets of $P$,
repeatedly applying the path selection theorem and the geometric reduction theorem.
At each step, we reduce the excess $S(X_i)-\chi(P)$
by a fixed proportion while controlling the loss in logarithmic size.

Put $\chi=\chi(P)$ and let $X_0=P$. For each nonempty set $X_i$ obtained during the iteration, put
\[
e_i=S(X_i)-\chi.
\]
We stop if
\begin{equation}\label{eq:sh01}
\log|X_i|\leqslant\chi
\qquad\text{or}\qquad
e_i\leqslant32\chi^{1/3}.
\end{equation}
Otherwise, let $d_i$ be the largest power of two not exceeding $e_i/16$. Then
\begin{equation}\label{eq:sh02}
\frac{e_i}{32}<d_i\leqslant\frac{e_i}{16},\qquad
d_i>\chi^{1/3}.
\end{equation}
By the ending-rank encoding, the points of $X_i$ are assigned distinct lattice paths
of length $S(X_i)\leqslant2\chi<2d_i^3$. Moreover, by~\eqref{eq:sh01} and~\eqref{eq:sh02},
\[
|X_i|>2^\chi>2^{S(X_i)-32d_i}.
\]
Theorem~\ref{thm:paths} therefore gives a nonempty subset $T_i\subseteq X_i$ whose associated paths have width at most $3d_i$ and such that
\[
\log\frac{|X_i|}{|T_i|}
=O\left(\sqrt{\frac{S(X_i)}{d_i}}+\log(d_i+1)\right).
\]
For $p,q\in T_i$ with $p\prec q$, the rank vertices corresponding to the slope $\sigma(p,q)$
are incomparable by~\eqref{eq:rank-increments}. Their height difference is $z_q(\sigma(p,q))-z_p(\sigma(p,q))$,
so the width bound gives
\[
z_q(\sigma(p,q))-z_p(\sigma(p,q))\leqslant3d_i,
\]
where the ranks are taken with respect to $X_i$.
Applying Theorem~\ref{prop:reduction}, we obtain a nonempty subset $X_{i+1}\subseteq T_i$ satisfying
\[
|X_{i+1}|\geqslant\frac{|T_i|}{3\lceil\chi(X_i)/d_i\rceil}
\]
and
\[
S(X_{i+1})\leqslant\frac{\chi(X_i)+S(X_i)}2+4d_i.
\]
Since $\chi(X_i)\leqslant\chi$, it follows from~\eqref{eq:sh02} that
\begin{equation}\label{eq:sh03}
e_{i+1}\leqslant\frac{e_i}{2}+4d_i\leqslant\frac34e_i.
\end{equation}
Thus the iteration terminates after finitely many steps.

We next bound the total loss in logarithmic size.
Combining the bounds from Theorems~\ref{thm:paths} and~\ref{prop:reduction}, we obtain
\begin{equation}\label{eq:sh04}
\begin{aligned}
\log\frac{|X_i|}{|X_{i+1}|}
& \leqslant O\left(\sqrt{\frac{S(X_i)}{d_i}}+\log(d_i+1)\right)+\log\bigl(3\lceil\chi/d_i\rceil\bigr)\\
& =O\left(\sqrt{\frac{\chi}{e_i}}+e_i^{1/3}\right),
\end{aligned}
\end{equation}
where the last equality follows from~\eqref{eq:sh02},
$S(X_i)\leqslant 2\chi$, and $0<e_i\leqslant\chi$. Indeed,
\[
\sqrt{\frac{S(X_i)}{d_i}}=O\left(\sqrt{\frac{\chi}{e_i}}\right),\qquad
\log(d_i+1)=O(e_i^{1/3}),
\]
and
\[
\log\bigl(3\lceil\chi/d_i\rceil\bigr)=O\bigl(\log(1+\chi/e_i)\bigr)=O\left(\sqrt{\frac{\chi}{e_i}}\right).
\]
The logarithmic estimates follow from the elementary bounds $\log(1+t)=O(t^{1/3})$ and $\log(1+t)=O(t^{1/2})$ for $t\geqslant1$.

Suppose that the iteration terminates at $X_k$, and put $\rho=3/4$.
If $k\geqslant1$, then repeated application of~\eqref{eq:sh03} gives
\[
\rho^{-(k-1-i)}e_{k-1}\leqslant e_i\leqslant\rho^ie_0\qquad(0\leqslant i<k).
\]
Since $e_0\leqslant\chi$ and $e_{k-1}>32\chi^{1/3}$, we have
\[
\begin{aligned}
\sum_{i=0}^{k-1}\sqrt{\frac{\chi}{e_i}} & \leqslant\frac{1}{1-\sqrt{\rho}}\sqrt{\frac{\chi}{e_{k-1}}}=O(\chi^{1/3}),\\
\sum_{i=0}^{k-1}e_i^{1/3} & \leqslant\frac{e_0^{1/3}}{1-\rho^{1/3}}=O(\chi^{1/3}).
\end{aligned}
\]
Consequently, by~\eqref{eq:sh04},
\begin{equation}\label{eq:sh05}
\sum_{i=0}^{k-1}\log\frac{|X_i|}{|X_{i+1}|}=O(\chi^{1/3}).
\end{equation}
This also holds when $k=0$.

By the stopping condition~\eqref{eq:sh01} and the encoding bound~\eqref{eq:encoding-bound}, we have
\[
\log|X_k|\leqslant\chi+32\chi^{1/3}.
\]
Together with~\eqref{eq:sh05}, this yields
\[
\log|P|=\log|X_k|+\sum_{i=0}^{k-1}\log\frac{|X_i|}{|X_{i+1}|}\leqslant\chi+O(\chi^{1/3}).
\]

Finally, if $P$ contains no $n$ points in convex position, then $\chi(P)\leqslant n-1$,
and so $|P|\leqslant2^{n-1+O(n^{1/3})}$. It follows that $\ES(n)\leqslant2^{n+O(n^{1/3})}$, completing the proof.
\end{proof}

\end{document}